\documentclass[]{amsart}
\usepackage{amssymb}
\usepackage{amsxtra}
\usepackage{amsmath}
\usepackage{amsfonts}

\def\R{{\mathbb R}}
\def\C{{\mathbb C}}

\def\diag{\mbox{diag \-}}

\begin{document}

\title[]{ $G_2$-Holonomy  Metrics
Connected with a~$3$-Sasakian Manifold }
 \author[]{Bazaikin ~Ya.~V.~ and ~Malkovich~E.~G. }

\address{Ya.~V. Bazaikin,
Sobolev Institute of Mathematics,
 Novosibirsk, Russia}

\address{E.~G. Malkovich,
Novosibirsk State University, Novosibirsk, Russia}

\email{bazaikin@math.nsc.ru, malkovicheugen@ngs.ru}

\date{}

\maketitle
\begin{abstract}We construct complete noncompact Riemannian
metrics with $G_2$-holonomy on noncompact orbifolds that are $\Bbb
R^3$-bundles with the twistor space $\mathcal{Z}$ as a~spherical
fiber.\\~

\end{abstract}
\section[]{Introduction}

This article addressing $G_2$-holonomy metrics  is a natural
continuation of the study of $Spin(7)$-holonomy metrics which was
started in \cite{Bazaikin}. We consider an arbitrary 7-dimensional
compact $3$-Sasakian manifold $M$ and discuss the existence of a
smooth resolution of the conic metric over the twistor space
$\mathcal{Z}$ associated with $M$.

Briefly speaking, a manifold $M$ is $3$-Sasakian if and only if
the standard metric on the cone over~ $M$ is hyper-Kahler. Each
manifold of this kind $M$ is closely related to the twistor space
$\mathcal{Z}$ which is an orbifold with a Kahler--Einstein metric.
We consider the metrics that are natural resolutions of the
standard conic metric over $\mathcal{Z}$:
$$
\bar{g}=dt^2+A(t)^2(\eta_2^2+\eta_3^2)+B(t)^2 (\eta_4^2+\eta_5^2)
+ C(t)^2 (\eta_6^2+\eta_7^2),\eqno{(*)}
$$
where $\eta_2$ and $\eta_3$ are the characteristic $1$-forms of
$M$, $\eta_4$, $\eta_5$, $\eta_6$, and $\eta_7$ are the forms that
annul the $3$-Sasakian foliation on $M$, and $A$, $B$, and $C$ are
real functions.

One of the main results of the article is the construction (in the
case when $M/SU(2)$ is Kahler) of a $G_2$-structure which is
parallel with respect to  $(\ast)$ if and only if the following
system of ordinary differential equations is satisfied:
$$
\begin{array}{l}
A'=\frac{2A^2-B^2-C^2}{BC},\\
B'=\frac{B^2-C^2-2A^2}{CA},\\
C'=\frac{C^2-2A^2-B^2}{AB}.
\end{array}\eqno{(**)}
$$
In case $(\ast\ast)$ we thus see that $(\ast)$ has holonomy $G_2$;
hence, $(\ast)$ is Ricci-flat. The system of equations $(\ast
\ast)$ was previously obtained in \cite{Gibbons} in the~particular
case $M=SU(3)/S^1$.

For a solution to $(**)$ to be defined on some orbifold or
manifold, some additional boundary conditions  are required at
$t_0$ that we will state them later. These conditions cannot be
satisfied unless $B=C$, which leads us to the functions that give
rise to the solutions found originally in \cite{Bryant-Solomon}
when $M=S^7$ and $M=SU(3)/S^1$. If $B=C$ then  $(*)$ is defined on
the total space of an $\R^3$-bundle $\mathcal{N}$ over a
quaternionic-Kahler orbifold $\mathcal{O}$. In general,
$\mathcal{N}$ is an orbifold except in the event that $M=S^7$ and
$M=SU(3)/S^1$. Note that it is unnecessary for $\mathcal{O}$ to be
Kahler  in case $B=C$.

Finally, we consider the well-known examples of the $3$-Sasakian
manifolds constructed in \cite{Boyer-Galicki-Mann} and describe
the topology of the corresponding orbifolds $\mathcal{N}$.

\section[]{Construction of a Parallel $G_2$-Structure}

The definition of $3$-Sasakian manifolds, their basic properties,
and further references can be found in \cite{Bazaikin}. We mainly
take our notation from \cite{Bazaikin}.

Let $M$ be a  7-dimensional compact $3$-Sasakian manifold with
characteristic fields $\xi^1$, $\xi^2$, and $\xi^3$ and
characteristic $1$-forms $\eta_1$, $\eta_2$, and $\eta_3$.
Consider the principal bundle $\pi: M \rightarrow \mathcal{O}$
with the structure group $Sp(1)$ or $SO(3)$ over the
quaternionic-Kahler orbifold $\mathcal{O}$ associated with $M$. We
are interested in the special case when $\mathcal{O}$ additionally
possesses a Kahler structure.

The field $\xi^1$ generates a locally free action of the circle
$S^1$ on $M$, and the metric on the twistor space
$\mathcal{Z}=M/S^1$ is a Kahler--Einstein metric. It is obvious
that $\mathcal{Z}$ is topologically a bundle over $\mathcal{O}$
with fiber $S^2=Sp(1)/S^1$ (or $S^2=SO(3)/S^1$) associated with
$\pi$. Consider the obvious action of $SO(3)$ on $\R^3$. The
two-fold cover $Sp(1) \rightarrow SO(3)$ determines the action of
$Sp(1)$ on $\R^3$, too. Now, let ${\mathcal{N}}$ be a bundle over
$\mathcal{O}$ with fiber $R^3$ associated with $\pi$. It is easy
to see that $\mathcal{O}$ is embedded in $\mathcal{N}$ as the zero
section, and $\mathcal{Z}$ is embedded in $\mathcal{N}$ as a
spherical section. The space $\mathcal{N} \backslash \mathcal{O}$
is diffeomorphic to the product $\mathcal{Z} \times (0,\infty)$.
Note that $\mathcal{N}$ can be assumed to be the projectivization
of the bundle $\mathcal{M}_1 \rightarrow {\mathcal{O}}$ of
\cite{Bazaikin}. In general, ${\mathcal{N}}$ is a 7-dimensional
orbifold; however, if $M$ is a regular $3$-Sasakian space then
$\mathcal{N}$ is a 7-dimensional manifold.

Let $\{ e^i \}, i=0,2,3, \dots, 7$, be an orthonormal basis of
$1$-forms on the standard Euclidean space ${\R}^7$ (the numeration
here is chosen so as to emphasize the connection with the
constructions of \cite{Bazaikin} and to keep the original notation
wherever possible). Putting $e^{ijk}=e^i \wedge e^j\wedge e^k$,
consider the following $3$-form $\Psi_0$ on ${\Bbb R}^7$:
$$
\Psi_0=-e^{023}-e^{045}+e^{067}+e^{346}-e^{375}-e^{247}+e^{256}.
$$

A differential $3$-form $\Psi$ on an oriented $7$-dimensional
Riemannian manifold $N$ defines a $G_2$-structure if, for each $p
\in N$, there exists an orientation-preserving isometry $\phi_p :
T_p N \rightarrow {\Bbb R}^7$ defined in a neighborhood of $p$
such that $\phi_p^*\Psi_0=\Psi|_p$. In this case the form $\Psi$
defines the unique metric $g_{\Psi}$ such that
$g_{\Psi}(v,w)=\langle \phi_p v, \phi_p w\rangle$ for $v,w \in T_p
N$ \cite{Bryant-Solomon}. If the form $\Psi$ is parallel $(\nabla
\Psi =0)$ then the holonomy group of the Riemannian manifold $N$
lies in $G_2$. The parallelness of the form $\Psi$ is equivalent
to its closeness and cocloseness \cite{Gray}:
$$
d \Psi=0,\quad  d \ast \Psi =0.
            \eqno{(1)}
$$

Note that the form $\Phi_0=e^1\wedge \Psi_0 - \ast \Psi_0$, where
$\ast$ is the Hodge operator in ${\R}^7$, determines a
$Spin(7)$-structure on ${\R}^8$ with the orthonormal basis
$\{e^i\}_{i=0,1,2,\dots ,7}$.

Locally choose an orthonormal system $\eta_4,\eta_5,\eta_6,\eta_7$
that generates the annihilator of the vertical subbundle
${\mathcal{V}}$ so that
$$
\omega_1=2( \eta_4 \wedge \eta_5 - \eta_6 \wedge \eta_7), \quad
\omega_2=2( \eta_4 \wedge \eta_6 - \eta_7 \wedge \eta_5), \quad
\omega_3=2( \eta_4 \wedge \eta_7 - \eta_5 \wedge \eta_6),
$$
where the forms $\omega_i$ correspond to the quaternionic-Kahler
structure on ${\mathcal{O}}$. It is clear that
$\eta_2,\eta_3,\dots,\eta_7$ is an orthonormal basis for $M$
annulling the one-dimensional foliation generated by $\xi^1$;
therefore, we can consider the metric of the following form on
$(0,\infty) \times {\mathcal{Z}}$:
$$
\bar{g}=dt^2+A(t)^2 \bigl(\eta_2^2+\eta_3^2\bigr) +B(t)^2
\bigl(\eta_4^2+\eta_5^2\bigr) + C(t)^2
\bigl(\eta_6^2+\eta_7^2\bigr).
              \eqno{(2)}
$$
Here  $A(t)$, $B(t)$, and $C(t)$ are defined on the interval
$(0,\infty)$.

We suppose that $\mathcal{O}$ is a Kahler orbifold; therefore,
$\mathcal{O}$ has the closed Kahler form that can be lifted to the
horizontal subbundle $\mathcal{H}$ as a closed form $\omega$.
Without loss of generality we can assume that we locally have
$$
\omega=2( \eta_4 \wedge \eta_5 + \eta_6 \wedge \eta_7).
$$
If we now put
$$
e^0=dt, \quad e^i=A \eta_i,\ i=2,3, \quad  e^j=B \eta_j,\ j=4,5,
\quad e^k=C \eta_k,\ k=6,7,
$$
then the forms $\Psi_0$ and $\ast \Psi_0$ become
$$
\Psi_1=-e^{023}-\frac{B^2+C^2}{4}e^0\wedge\omega_1
-\frac{B^2-C^2}{4}e^0\wedge\omega+\frac{BC}{2}e^3\wedge\omega_2-
\frac{BC}{2}e^2\wedge\omega_3,
$$
$$
\Psi_2=C^2B^2\Omega-\frac{B^2+C^2}{4}e^{23}\wedge\omega_1-\frac{B^2-C^2}{4}e^{23}\wedge\omega+\frac{BC}{2}e^{02}\wedge\omega_2
+\frac{BC}{2}e^{03}\wedge\omega_3,
$$
where $\Omega=\eta_4 \wedge \eta_5 \wedge \eta_6 \wedge
\eta_7=-\frac{1}{8} \omega_1 \wedge \omega_1=-\frac{1}{8} \omega_2
\wedge \omega_2=-\frac{1}{8} \omega_3 \wedge \omega_3$.

It is now obvious that $\Psi_1$ and $\Psi_2$ are defined globally
and independently of the local choice of $\eta_i$; consequently,
they  uniquely define the metric $\bar{g}$ given locally by $(2)$.
Then the condition $(1)$ that the holonomy group lies in $G_2$ is
equivalent to the equation
$$
d\Psi_1=d\Psi_2=0. \eqno{(3)}
$$

{\bf Theorem}.{\it If $\mathcal{O}$ possesses a Kahler structure
then $(2)$ on $\mathcal{N}$ is a smooth metric with holonomy $G_2$
given by the form $\Psi_1$ if and only if the functions $A$, $B$,
and $C$ defined on the interval $[t_0,\infty)$ satisfy the system
of ordinary differential equations
$$
A'=\frac{2A^2-B^2-C^2}{BC}, \quad B'=\frac{B^2-C^2-2A^2}{CA},
\quad C'=\frac{C^2-2A^2-B^2}{AB}
              \eqno{(4)}
$$
with the initial conditions

1. \quad $A(0)=0$ and $|A'(0)|= 2$;

2. \quad $B(0), C(0) \neq 0$, and $B'(0)=C'(0)=0$;

3. \quad the functions $A$, $B$, and $C$ have fixed sign on the
interval $(t_0,\infty)$. }

\vskip0.2cm

{\bf Proof}.

In \cite{Bazaikin} the following relations were obtained, closing
 the algebra of forms:
$$
\begin{array}{l} d e^0 =0,
\\
d e^i =\frac{A_i'}{A_i} e^0 \wedge e^i+ A_i \omega_i- \frac{2
A_i}{A_{i+1} A_{i+2}} e^{i+1} \wedge e^{i+2}, \quad i=1,2,3\quad
\textrm{mod} \quad 3,
\\
d \omega_i  = \frac{2}{A_{i+2}} \omega_{i+1} \wedge e^{i+2}
-\frac{2}{A_{i+1}} e^{i+1} \wedge \omega_{i+2},  \quad i=1,2,3
\quad \textrm{mod} \quad 3.
\end{array}
$$
By adding the relation $d \omega=0$ and carrying out some
calculations to be omitted here, we obtain the sought system.

The smoothness conditions for the metric at  $t_0$ are proven by
analogy with the case of holonomy $Spin(7)$ which was elaborated
in \cite{Bazaikin}. We only note that, taking the quotient of the
unit sphere $S^3$ by the Hopf action of the circle, we obtain the
sphere of radius $1/2$, which explains the condition $|A'(0)|=2$.

In case $B=C$ the system reduces to the pair of equations
$$
A'=2 \left( \frac{A^2}{B^2}-1 \right), \quad B'=-2\frac{A}{B}
$$
whose solution gives the metric
$$
\bar{g}=\frac{dr^2}{1-{r_0^4}/{r^4}}+r^2 \left(
1-\frac{r_0^4}{r^4} \right) \bigl( \eta_2^2+\eta_3^2 \bigr)+2 r^2
\bigl( \eta_4^2+\eta_5^2 + \eta_6^2+\eta_7^2 \bigr).
$$
The regularity conditions hold. This smooth metric was originally
found in \cite{Bryant-Solomon} in the event that $M=SU(3)/S^1$ and
$M=S^7$ (observe that we need not require $\mathcal{O}$ to be
Kahler when $B=C$).

In the general case $B \neq C$ system $(4)$ can also be integrated
\cite{Gibbons}. However, the resulting solutions do not enjoy the
regularity conditions.

\section[]{Examples}

Some interesting family of examples arises when we consider the
7-dimensional biquotients of the Lie group $SU(3)$ as $3$-Sasakian
manifolds. Namely, let $p_1$, $p_2$, and $p_3$ be pairwise coprime
positive integers. Consider the following action of $S^1$ on the
Lie group $SU(3)$:
$$
z \in S^1: A \mapsto \diag(z^{p_1}, z^{p_2}, z^{p_3}) \cdot A
\cdot \diag(1,1,z^{-p_1-p_2-p_3}).
$$
This action is free; moreover, it was demonstrated in
\cite{Boyer-Galicki-Mann} that there is a $3$-Sasakian structure
on the orbit space ${\mathcal{S}}={\mathcal{S}}_{p_1,p_2,p_3}$.
Moreover, the action of $SU(2)$ on $SU(3)$ by right translations

$$
B \in SU(2): A \mapsto A \cdot \left( \begin{array}{cc} B & 0 \\ 0
& 1 \end{array} \right)
$$

commutes with the action of $S^1$ and can be pushed forward to the
orbit space $\mathcal{S}$. The corresponding Killing fields will
be the characteristic fields $\xi_i$ on $\mathcal{S}$. Therefore,
the corresponding twistor space
$\mathcal{Z}={\mathcal{Z}}_{p_1,p_2,p_3}$ is the orbit space of
the following action of the torus $T^2$ on $SU(3)$:
$$
(z,u) \in T^2: A \mapsto \diag(z^{p_1}, z^{p_2}, z^{p_3}) \cdot A
\cdot \diag(u,u^{-1},z^{-p_1-p_2-p_3}).
              \eqno{(5)}
$$

{\bf Lemma}.{\it The space ${\mathcal{Z}}_{p_1,p_2,p_3}$ is
diffeomorphic to the orbit space of  $U(3)$ with respect to the
following action of $T^3$:
$$
(z,u,v) \in T^3: A \mapsto
\diag(z^{-p_2-p_3},z^{-p_1-p_3},z^{-p_1-p_2}) \cdot A \cdot
\diag(u,v,1).
              \eqno{(6)}
$$   }

It suffices to verify that each $T^3$-orbit in $U(3)$ exactly cuts
out an orbit of the $T^2$-action $(5)$ in $SU(3) \subset U(3)$.

Action $(6)$ makes it possible to describe the topology of
$\mathcal{Z}$ and, consequently, the topology of $\mathcal{N}$
clearly. Here we use the construction of \cite{Eschenburg}.
Consider the submanifold $E=\{ (u,[v]) \mid u \bot v \} \subset
S^5 \times {\C} P^2.$ It is obvious that $E$ is diffeomorphic to
$U(3)/S^1\times S^1$ (the ''right'' part of $(6)$) and is the
projectivization of the ${\C}^2$-bundle $\widetilde{E}=\{
(u,v)\mid  u \bot v \} \subset S^5 \times {\C}^3$ over $S^5$. By
adding the trivial one-dimensional complex bundle over $S^5$ to
$\widetilde{E}$, we obtain the trivial bundle $S^5 \times {\C}^3$
over $S^5$.

The group $S^1$ acts from the left by the automorphisms of the
vector bundle $\widetilde{E}$, and $\mathcal{Z}=S^1 \backslash E$
is the projectivization of the ${\C}^2$-bundle $S^1 \backslash
\widetilde{E}$ over the weighted complex projective space
$\mathcal{O}={\C} P^2 (q_1,q_2,q_3)=S^1 \backslash S^5$, where
$q_i=(p_{i+1}+p_{i+2})/2$ for $p_i$ all odd and
$q_i=(p_{i+1}+p_{i+2})$ otherwise.

The above implies that the bundle $S^1 \backslash \widetilde{E}$
is stably equivalent to the bundle $S^1 \backslash (S^5 \times
{\C}^3)$ over ${\mathcal{O}}$. The last bundle splits obviously
into the Whitney sum $\sum\nolimits_{i=1}^3 \xi^{q_i}$, where
$\xi$ is an analog of the one-dimensional universal bundle of
$\mathcal{O}$.

{\bf Corollary}.{\it The twistor space $\mathcal{Z}$ is
diffeomorphic to the projectivization of a two-dimensional complex
bundle over $\C P^2(q_1,q_2,q_3)$ which is stably equivalent to
$\xi^{q_1} \oplus \xi^{q_2} \oplus \xi^{q_3}$. }

\end{document}